\documentclass[11pt]{article}

\usepackage{amsmath}
\usepackage{times,theorem,latexsym,amssymb}
\usepackage[multiple]{footmisc}
\newcommand{\pagestart}{1}

\makeatletter

\def\@evenhead{\footnotesize\thepage\hfil\slshape\leftmark}%
\def\@oddhead{\footnotesize{\slshape\rightmark}\hfil\thepage}%

\numberwithin{equation}{section} \numberwithin{figure}{section}
\numberwithin{table}{section}

\renewcommand{\thefootnote}{\roman{footnote}}
\renewcommand{\thefootnote}{\Roman{footnote}}
\renewcommand{\thefootnote}{\arabic{footnote}}

\makeatother

\newenvironment{summary}{\vskip\baselineskip \noindent\small\bf Summary: \rm}%
{\vskip\baselineskip}

\newtheorem{theorem}{Theorem }[section]

\newtheorem{definition}[theorem]{Definition}
{\theorembodyfont{\rmfamily}}
\newtheorem{lemma}[theorem]{Lemma}

{\theorembodyfont{\rmfamily}}

\title{Parametric Summability and Its Applications}
\author{Jinlu Li, Robert Mendris}
\date{April 11, 2015}

\newcommand{\hide}[1]{}

\newcommand{\br}{\newline\noindent}

\newcommand{\mylabel}[1]{\label{#1}}  



\AtBeginDocument{

}

\begin{document}
\maketitle\thispagestyle{empty}


\begin{summary}
{ In this paper we study summability based on double sequences of
complex constants as it is defined in {}``Linear Operators, General
Theory'' by N. Dunford and J. T. Schwartz. We define {}``power double
sequences'' or infinite {}``power matrices'' as certain generalizations
of double sequences and power series. }

{ We relate the summability and boundedness of the power double sequences
to the summability and boundedness of the double sequence \textit{A.}
}

{ While others do investigate {}``power matrices'' their definitions,
as far as we were able to find, differ from our definitions. Using
these definitions we extend some summability results for double sequences
of constants to our power double sequences.}

{We investigate also other possible generalizations of double sequences
and assess their usefulness in summability study.} 

\end{summary}

\renewcommand{\thefootnote}{}
\footnotetext{\hspace*{-.51cm}
}

\section{Introduction}\mylabel{sec:1}

The only (most recent) related work we found is ...

There are different reasons for transforming one sequence into another and
  most of them are related to convergence.
  The practical need to improve convergence gave the impulz to study sequence
  transformations already in 17th century and resulted in the creation of
  summability theory at the end of 19th century.
  Before the invention of computers, mainly linear sequence transformations
  were studied.
  Approaches based on classical analysis culminated when [1] was published. 
  After that modern approaches based on functional analysis appeared.
  For a comprehensive review of classical and modern methods in summability
  see [2].
  From practical point of view, regular linear transformations are in general
  at most moderately powerful in improving convergence, and the popularity
  of most linear transformations has declined considerably in recent years.
  It seems, however, that the limiting factor is regularity not linearity.
  Recently also new powerful non-linear sequence transformations attracted
  research and applications.
  This is discussed in a nice historical review [3].

 In this paper, we contribute to the classical summability methods of double
  sequences.   
  As a review of these methods see [4].

{ Let $A = \{ a_{{ij}} \},\; i = 1, 2, \ldots{},\; j = 0, 1, 2, \ldots{}$,
be a double sequence of complex constants, that is, 
}

\begin{center}
{ $A={\left(\begin{matrix}a_{10} & a_{11} & a_{12} & \text{......}\\
a_{20} & a_{21} & a_{22} & \text{......}\\
a_{30} & a_{31} & a_{32} & \text{......}\\
\text{...} & \text{...} & \text{...} & \text{......}\end{matrix}\right)}$ } 
\par\end{center}

{ Let $\Delta$ be the set of all double sequences of complex constants
and $c$ be the space of all convergent sequences of scalars.
}

\begin{definition}\mylabel{def:regsummability}
  Suppose that a double sequence $\{a_{ij}\}$
defines a linear transformation $T$ of $c$ onto itself by means of formula 
\[
T[s_{0},s_{1},\ldots]=[t_{0},t_{1},\ldots]=[\overset{{\infty}}{\underset{{i=0}}{\sum}}a_{0i}s_{i},\overset{{\infty}}{\underset{{i=0}}{\sum}}a_{1i}s_{i},\;\ldots].
\]
If $T$ preserves limits of sequences 
(i.e. if $\underset{{i\rightarrow\infty}}{{\text{lim}}}t_{{i}}=\underset{{i\rightarrow\infty}}{{\text{lim}}}s_{{i}}$
for every $[s_{i}]\in c$), then the double sequence (or matrix) $A=\{a_{ij}\}$
is said to define a \emph{regular method of summability}.
\end{definition}

\begin{lemma}\mylabel{lem:1}
\textit{A}(\textit{a}) defines a bounded linear
map of \textit{c} into \textit{c}, if and only if the following three
conditions hold:
\begin{enumerate}
\item { ${\underset{{1\le i<\infty}}{{\text{sub}}}\overset{{\infty}}{\underset{{j=0}}{\sum}}{|f_{{ij}}(a)|}=M<\infty}$;} 
\item {${\underset{{i\rightarrow\infty}}{{\text{lim}}}f_{{ij}}(a)}$ exists
for \textit{j} = 1, 2, \ldots{}. ;} 
\item { ${\underset{{i\rightarrow\infty}}{{\text{lim}}}\overset{{\infty}}{\underset{{j=0}}{\sum}}{f_{{ij}}(a)}}$
exists.} 
\end{enumerate}
 \end{lemma}
 
{ Once \textit{A}(\textit{a}) satisfies the above three conditions,
for any column vector $s=[{s_{{1}},s_{{2}},s_{{3}},\;\ldots}]$ $\in c$,
the bounded linear map, \textit{A}(\textit{a}), of \textit{c} into
\textit{c} is defined by}
\[
{A(a){(s)_{i}}={f_{i0}(a)\underset{i\rightarrow\infty}{\text{lim}}s_{i}+\overset{\infty}{\underset{j=1}{\sum}}{f_{ij}(a)s_{j}}}\textrm{, for }i=1,2,\ldots,}
\]

{ where \textit{A}(\textit{a})${(s)}$ = {[}\textit{A}(\textit{a})${(s)_{{1}}}$,
\textit{A}(\textit{a})${(s)_{{2}}}$, \textit{A}(\textit{a})${(s)_{{3}}}$,
\ldots{}]. }

\bigskip{}

{ Let ${B_{{c}}}$ denote the space of all linear bounded maps of
\textit{c} into \textit{c}. From Lemma 1, if \textit{A}(\textit{a})
$\in{B_{{c}}}$ , then it has norm }

\[  {|A(a)|=M}  \]

\begin{theorem}\mylabel{th:Toeplitz}
{ \textbf{Silverman-Toeplitz.} \textit{A}(\textit{a}) defines
a regular method of summability, if and only if the following three
conditions hold:}
\begin{enumerate}
\item { ${\underset{1\le i<\infty}{\text{lub}}\overset{\infty}{\underset{j=0}{\sum}}{|f_{ij}(a)|}=M<\infty}$;} 
\item { ${\underset{i\rightarrow\infty}{\text{lim}}f_{ij}(a)}= 0$ for $j = 1, 2, \ldots{}$;} 
\item { ${\underset{i\rightarrow\infty}{\text{lim}}\overset{\infty}{\underset{j=1}{\sum}}{f_{ij}(a)}}= 1.$} 
\end{enumerate}
\end{theorem}

{ It is clear that if \textit{A}(\textit{a}) defines a regular method
of summability, then \textit{A}(\textit{a}) $\in{B_{c}}$.}

One can find all the Definition, the Lemma 1. and Silverman-Toeplitz
Theorem in {[}Dunford-Schwartz]. ({}``Linear Operators, General Theory''
by N. Dunford and J. T. Schwartz.) 

\bigskip{}

{ Let $A(z) = \{ f_{ij}(z) \},\; i= 1, 2, \ldots{},\; j = 0, 1, 2, \ldots{}$,
be a double sequence of functions with same domain $D$, 
which is a subset of complex numbers, that is, 
}

\begin{center}
{$A(z)={\left(\begin{matrix}f_{10}(z) & f_{11}(z) & f_{12}(z) & \text{......}\\
f_{20}(z) & f_{21}(z) & f_{22}(z) & \text{......}\\
f_{30}(z) & f_{31}(z) & f_{32}(z) & \text{......}\\
\text{...} & \text{...} & \text{...} & \text{......}\end{matrix}\right)}$ .} 
\par\end{center}

{ It is clear that for any input ${a\in D}$, the output of $A(a)\in\Delta$.
}

\section{Preliminary and Review}\mylabel{sec:2}

\subsection*{Power Matrices}\mylabel{sec:3}

Let \textit{A} = \{ ${a_{{ij}}}$ \}, \textit{i} = 1, 2, \ldots{}
, \textit{j} = 0, 1, 2, \ldots{} , be a double sequence of complex
constants.

{ The column power matrix induced by \textit{A} is defined as}

\[
{P_{A}^{C}(z)=\{{a_{{ij}}z^{{i}}}\},i=1,2,\ldots,j=0,1,2,\ldots\textrm{, that is,}}\]

\bigskip{}

\begin{center}
{${P_{{A}}^{{C}}(z)}={\left(\begin{matrix}a_{10}z & a_{11}z & a_{12}z & \text{......}\\
a_{20}z^{{2}} & a_{21}z^{{2}} & a_{22}z^{{2}} & \text{......}\\
a_{30}z^{{3}} & a_{31}z^{{3}} & a_{32}z^{{3}} & \text{......}\\
\text{...} & \text{...} & \text{...} & \text{......}\end{matrix}\right)}$ .} 
\par\end{center}

\bigskip{}

{ The row power matrix induced by \textit{A} is defined as}

\[
{{P_{{A}}^{{R}}(z)}=\{{a_{{ij}}z^{{j}}}\},i=1,2,\ldots,j=0,1,2,\ldots\textrm{, that is,}}\]

\bigskip{}

\begin{center}
{${P_{{A}}^{{R}}(z)}={\left(\begin{matrix}a_{10} & a_{11}z & a_{12}z^{{2}} & \text{......}\\
a_{20} & a_{21}z & a_{22}z^{{2}} & \text{......}\\
a_{30} & a_{31}z & a_{32}z^{{2}} & \text{......}\\
\text{...} & \text{...} & \text{...} & \text{......}\end{matrix}\right)}$ .} \bigskip{}

\par\end{center}

{The double power function matrix induced from \textit{A} is defined
as}

\[
{{P_{{A}}(z)}=\{{a_{{ij}}z^{{i+j}}}\},i=1,2,\ldots,j=0,1,2,\ldots\textrm{, that is,}}\]
\bigskip{}

\begin{center}
{${P_{{A}}(z)}={\left(\begin{matrix}a_{10}z & a_{11}z^{{2}} & a_{12}z^{{3}} & \text{......}\\
a_{20}z^{{2}} & a_{21}z^{{3}} & a_{22}z^{{4}} & \text{......}\\
a_{30}z^{{3}} & a_{31}z^{{4}} & a_{32}z^{{5}} & \text{......}\\
\text{...} & \text{...} & \text{...} & \text{......}\end{matrix}\right)}$ .} \bigskip{}

\par\end{center}

We will immediately generalize these definitions in the following
section.\bigskip{}

\subsection*{General Power Matrices}\mylabel{sec:4}

We will define power matrices of the first type now.

{ \textbf{Definition.} Let $A=\{a_{ij}\}$, $i=0,1,2,\ldots$, $j=0,1,2,\ldots$,
be a double sequence of complex constants. Let ${g(z)=\overset{{\infty}}{\underset{{i=0}}{\sum}}{g_{{i}}z^{i}}}$
be a complex power series. Denote its radius of convergence as $R_{g}$.
The \emph{column power matrix} induced by \textit{A} and associated
with ${g}(z)$ is defined as}

{ ${P_{{A;g}}^{{C}}(z)}=\{{a_{{ij}}g_{{i}}z^{{i}}}\}$, $i=0,1,2,\ldots$,
$j=0,1,2,\ldots$ , that is, }

\begin{center}
{$P_{{A;g}}^{{C}}(z)=\left(\begin{array}{cccc}
a_{00}g_{0} & a_{01}g_{0} & a_{02}g_{0} & \ldots\\
a_{10}g_{1}z & a_{11}g_{1}z & a_{12}g_{1}z & \ldots\\
a_{20}g_{2}z^{2} & a_{21}g_{2}z^{2} & a_{22}g_{2}z^{2} & \ldots\\
\ldots & \ldots & \ldots & \ldots\end{array}\right)$ .} 
\par\end{center}

\bigskip{}

{ \textbf{Definition.} Let $A=\{a_{ij}\}$, $i=0,1,2,\ldots$, $j=0,1,2,\ldots$,
be a double sequence of complex constants. Let ${h(z)=\overset{{\infty}}{\underset{{j=0}}{\sum}}{h_{{j}}z^{j}}}$
be a complex power series. Denote its radius of convergence as $R_{h}$.
The \emph{row power matrix} induced by \textit{A} and associated with
${h}(z)$ is defined as}

{ ${P_{{A;h}}^{{R}}(z)}=\{{a_{{ij}}h_{{j}}z^{{j}}}\}$, $i=0,1,2,\ldots$,
$j=0,1,2,\ldots$ :, that is, }

\begin{center}
$P_{{A;h}}^{{R}}(z)=\left(\begin{array}{cccc}
a_{00}h_{0} & a_{01}h_{1}z & a_{02}h_{2}z^{2} & \ldots\\
a_{10}h_{0} & a_{11}h_{1}z & a_{12}h_{2}z^{2} & \ldots\\
a_{20}h_{0} & a_{21}h_{1}z & a_{22}h_{2}z^{2} & \ldots\\
\ldots & \ldots & \ldots & \ldots\end{array}\right)$ 
\par\end{center}

\bigskip{}

\noindent We will generalize the double power function matrix ${P_{{A}}(z)}$
now.

{ \textbf{Definition.} Let $A=\{a_{ij}\}$, $i=0,1,2,\ldots$, $j=0,1,2,\ldots$,
be a double sequence of complex constants. Let ${g(z)=\overset{{\infty}}{\underset{{i=0}}{\sum}}{g_{{i}}z^{i}}}$
be a complex power series. Denote its radius of convergence respectively
as $r_{g}$. The \textit{power double sequence} of \emph{second} type
induced by \textit{A} and associated with ${g}(z)$ is defined as}

{ ${P_{{A;g}}(z)}=\{{a_{{ij}}g_{{i+j}}z^{{i+j}}}\}$, $i=0,1,2,\ldots$,
$j=0,1,2,\ldots$ :}

\begin{center}
$P_{{A;g}}(z)=\left(\begin{array}{cccc}
a_{00}g_{0} & a_{01}g_{1}z & a_{02}g_{2}z^{2} & \ldots\\
a_{10}g_{1}z & a_{11}g_{2}z^{2} & a_{12}g_{3}z^{3} & \ldots\\
a_{20}g_{2}z^{2} & a_{21}g_{3}z^{3} & a_{22}g_{4}z^{4} & \ldots\\
\ldots & \ldots & \ldots & \ldots\end{array}\right)$ 
\par\end{center}

\bigskip{}
{ \textbf{Definition.} Let $A=\{a_{ij}\}$, $i=0,1,2,\ldots$, $j=0,1,2,\ldots$,
be a double sequence of complex constants. Let ${g(z)=\overset{{\infty}}{\underset{{i=0}}{\sum}}{g_{{i}}z^{i}}}$
and ${h(z)=\overset{{\infty}}{\underset{{j=0}}{\sum}}{h_{{j}}z^{j}}}$
be two complex power series. Denote their radius of convergence respectively
as $r_{g}$ and $r_{h}$. The \textit{power double sequence} of \emph{third}
type induced by \textit{A} and associated with ${g}(z)$ and ${h}(z)$
is defined as}

{ ${P_{{A;g,h}}(z)}=\{{a_{{ij}}g_{{i}}h_{{j}}z^{{i+j}}}\}$, $i=0,1,2,\ldots$,
$j=0,1,2,\ldots$ :}

\begin{center}
$P_{{A;g,h}}(z)=\left(\begin{array}{cccc}
a_{00}g_{0}h_{0} & a_{01}g_{0}h_{1}z & a_{02}g_{0}h_{2}z^{2} & \ldots\\
a_{10}g_{1}h_{0}z & a_{11}g_{1}h_{1}z^{2} & a_{12}g_{1}h_{2}z^{3} & \ldots\\
a_{20}g_{2}h_{0}z^{2} & a_{21}g_{2}h_{1}z^{3} & a_{22}g_{2}h_{2}z^{4} & \ldots\\
\ldots & \ldots & \ldots & \ldots\end{array}\right)$ 
\par\end{center}

\bigskip{}
\textbf{Remark}. More general definition would consider $\{g_{i}\}$
and $\{h_{j}\}$ to be two arbitrary number sequences.\bigskip{}

\section{Summability Results}\mylabel{sec:5}

For power double sequences of first type we have:

\bigskip{}
{ \textbf{Proposition 1}(C). Let the double sequence of complex constants
$\{a_{ij}\}$, $i=0,1,2,\ldots$, $j=0,1,2,\ldots$, be a regular
method of summability and ${g(z)=\overset{{\infty}}{\underset{{i=0}}{\sum}}{g_{{i}}z^{i}}}$
be a complex function with its power series. Then the following two
conditions are equivalent for any complex number $z$:

(\emph{i}) $\underset{{i\rightarrow\infty}}{{\text{lim}}}g_{{i}}z^{{i}}=1$.

(\emph{ii}) The power double sequence of first type $\left\{ \left(P_{{A;g}}^{{C}}(z)\right)_{ij}\right\} $
is a regular method of summability.

\textbf{Proof.} (\emph{i}) implies (\emph{ii}) is a straightforward
verification of the three conditions of Silverman-Toeplitz Theorem.

\noindent (\emph{ii}) implies (\emph{i}). The third condition of Silverman-Toeplitz
Theorem for $P_{{A;g}}^{{C}}$ and for $A$ gives (\emph{i}).\bigskip{}

\textbf{Proposition 1}(R). Let the double sequence of complex constants
$\{a_{ij}\}$, $i=0,1,2,\ldots$, $j=0,1,2,\ldots$, be a regular
method of summability and ${h(z)=\overset{{\infty}}{\underset{{j=0}}{\sum}}{h_{{j}}z^{j}}}$
be a complex function with its power series. Then the following two
conditions are equivalent for any complex number $z$, which the sequence
$\{h_{{i}}z^{{i}}\}$ is convergent for:

(\emph{i}) $\underset{{j\rightarrow\infty}}{{\text{lim}}}h_{{j}}z^{{j}}=1$.

(\emph{ii}) The power double sequence of first type $\left\{ \left(P_{{A;h}}^{{R}}(z)\right)_{ij}\right\} $
is a regular method of summability.

\textbf{Proof.} (\emph{i}) implies (\emph{ii}) is a straightforward
verification of the three conditions of Silverman-Toeplitz Theorem.

\noindent (\emph{ii}) implies (\emph{i}). The third condition of Silverman-Toeplitz
Theorem for $P_{{A;h}}^{{R}}$ and for $A$ along with the convergence
of the sequence $\{h_{{i}}z^{{i}}\}$ gives (\emph{i}) (see the proof
of (\emph{ii}) implies (\emph{i}) in the Proposition 1(II)).

\noindent \bigskip{}

\noindent For power double sequences of second type we have:

\bigskip{}

\textbf{Theorem} 1(II). Let double sequence of complex constants $\{a_{ij}\}$,
$i=0,1,2,\ldots$, $j=0,1,2,\ldots$, be a regular method of summability
and ${g(z)=\overset{{\infty}}{\underset{{i=0}}{\sum}}{g_{{i}}z^{i}}}$
be a complex function with its power series. Then the following two
conditions are equivalent for any complex number $z$, which the sequence
$\{g_{{i}}z^{{i}}\}$ is convergent for:

(\emph{i}) $\underset{{i\rightarrow\infty}}{{\text{lim}}}g_{{i}}z^{{i}}=1$.

(\emph{ii}) The power double sequence of first type $\left\{ \left(P_{{A;g}}(z)\right)_{ij}\right\} $
is a regular method of summability.

\bigskip{}

\textbf{Proof}. Let's show first that (\emph{ii}) implies (\emph{i}).
From (\emph{ii}) and the condition 3. of Silverman-Toeplitz theorem,
we have

\noindent ${\underset{{i\rightarrow\infty}}{{\text{lim}}}\overset{{\infty}}{\underset{{j=0}}{\sum}}{a_{{ij}}g_{{i+j}}z^{{i+j}}}}=1$.
Set $k=i+j$. This changes into $\underset{{i\rightarrow\infty}}{{\text{lim}}}\overset{{\infty}}{\underset{{k=i}}{\sum}}{a_{{i,k-i}}g_{{k}}z^{{k}}}=1$.
Set $b_{ik}=a_{i,k-i}$ for $i\leq k$ and zero otherwise. Observe
$\{b_{ik}\}$ is also a regular method of summability. Then by Silverman-Toeplitz
theorem the limits of convergent sequences are preserved: $\underset{{k\rightarrow\infty}}{{\text{lim}}}g_{{k}}z^{{k}}=\underset{{i\rightarrow\infty}}{{\text{lim}}}\overset{{\infty}}{\underset{{k=0}}{\sum}}{b_{{i,k}}g_{{k}}z^{{k}}}$.
But the right hand side is equal to $\underset{{i\rightarrow\infty}}{{\text{lim}}}\overset{{\infty}}{\underset{{k=i}}{\sum}}{a_{{i,k-i}}g_{{k}}z^{{k}}}=1$
and (\emph{i}) immediately follows.

Now we show that (\emph{i}) implies (\emph{ii}). We will use Silverman-Toeplitz
theorem again and we need to prove its three conditions:

1. ${\underset{{0\le i<\infty}}{{\text{sup}}}\overset{{\infty}}{\underset{{j=0}}{\sum}}{|\left(P_{{A;g}}(z)\right)_{ij}|}=\underset{{0\le i<\infty}}{{\text{sup}}}\overset{{\infty}}{\underset{{j=0}}{\sum}}{|a_{{ij}}g_{{i+j}}z^{{i+j}}|}}$
by the definition of $P_{{A;g}}$.

\noindent Set $k=i+j$. Then the above supremum is
  \\ $\underset{{0\le i<\infty}}{{\text{sup}}}\overset{{\infty}}{\underset{{k=i}}{\sum}}{|a_{{i,k-i}}g_{{k}}z^{{k}}|}\leq\underset{{0\le i<\infty}}{{\text{sup}}}\overset{{\infty}}{\underset{{k=i}}{\sum}}{|a_{{i,k-i}}|}\cdot\underset{{0\le k<\infty}}{{\text{sup}}}{|g_{{k}}z^{{k}}|}<\infty$
  \\ since the first supremum is finite by the first condition of $\{a_{ij}\}$
being a regular method of summability and the second one by the existence
of the limit in (\emph{i}).

\bigskip{}

2. ${\underset{{i\rightarrow\infty}}{{\text{lim}}}\left(P_{{A;g}}(z)\right)_{ij}}={\underset{{i\rightarrow\infty}}{{\text{lim}}}a_{{ij}}g_{{i+j}}z^{{i+j}}}$
for $j=0,1,2,\ldots$ by the definition.

\noindent Set $k=i+j$. Then the absolute value of the above limit
is $\left|\underset{{k\rightarrow\infty}}{{\text{lim}}}a_{{k-j,j}}g_{{k}}z^{{k}}\right|\leq\underset{{k\rightarrow\infty}}{{\text{lim}}}\left|a_{{k-j,j}}\right|\cdot\underset{{k\rightarrow\infty}}{{\text{lim}}}\left|g_{{k}}z^{{k}}\right|=0$
since the first limit is zero by the second condition of $\{a_{ij}\}$
being a regular method of summability and the second limit is $1$
by (\emph{i}).

\bigskip{}

3. ${\underset{{i\rightarrow\infty}}{{\text{lim}}}\overset{{\infty}}{\underset{{j=0}}{\sum}}{\left(P_{{A;g}}(z)\right)_{ij}}}={\underset{{i\rightarrow\infty}}{{\text{lim}}}\overset{{\infty}}{\underset{{j=0}}{\sum}}{a_{{ij}}g_{{i+j}}z^{{i+j}}}}$
by the definition.

\noindent Set $k=i+j$. Starting with (\emph{i}) following the fist
part of this proof ((\emph{ii}) implies (\emph{i})) backwards we have

\noindent $1=\underset{{k\rightarrow\infty}}{{\text{lim}}}g_{{k}}z^{{k}}=\underset{{i\rightarrow\infty}}{{\text{lim}}}\overset{{\infty}}{\underset{{k=0}}{\sum}}{b_{{i,k}}g_{{k}}z^{{k}}}=\underset{{i\rightarrow\infty}}{{\text{lim}}}\overset{{\infty}}{\underset{{k=i}}{\sum}}{a_{{i,k-i}}g_{{k}}z^{{k}}}={\underset{{i\rightarrow\infty}}{{\text{lim}}}\overset{{\infty}}{\underset{{j=0}}{\sum}}{a_{{ij}}g_{{i+j}}z^{{i+j}}}}$.

And this finishes the proof that $\left\{ \left(P_{{A;g}}(z)\right)_{ij}\right\} $
is a regular method of summability by Silverman-Toeplitz theorem.

\bigskip{}

\textbf{Remarks}. The condition (\emph{i}) in Theorem A. implies $|z|=r_{g}$,
where $r_{g}$ is the radius of convergence. The requirement in the
Theorem A. that the sequence $\{g_{{i}}z^{{i}}\}$ must be convergent
seems to be too restrictive but the condition (\emph{ii}) does not
guarantee its convergence. There are examples of non-convergent sequences
$\{g_{{i}}z^{{i}}\}$ (for both bounded and unbounded case) and regular
methods of summability that map these sequences to convergent ones
(= that sum them). Then by choosing $z=1$ one gets a counterexample
for each case.

\bigskip{}

\textbf{Corollary}. From the proof of the Theorem A. it is clear that
for $|z|<r_{g}$, conditions 1. and 2. hold but the limit in 3. is
zero and we don't get a regular method of summability in that case.

\bigskip{}

For power double sequences of third type we have:

$\vphantom{}$

\textbf{Theorem} 1(III). Let double sequence of complex constants
$\{a_{ij}\}$, $i=0,1,2,\ldots$, $j=0,1,2,\ldots$, be a regular
method of summability, and ${g(z)=\overset{{\infty}}{\underset{{i=0}}{\sum}}{g_{{i}}z^{i}}}$
and ${h(z)=\overset{{\infty}}{\underset{{j=0}}{\sum}}{h_{{j}}z^{j}}}$
be two complex power series. If $\underset{{i\rightarrow\infty}}{{\text{lim}}}g_{{i}}z^{{i}}$
and $\underset{{j\rightarrow\infty}}{{\text{lim}}}h_{{j}}z^{{j}}$
exist then the following two conditions are equivalent for any such
complex number $z$:

(\emph{i}) $\underset{{i\rightarrow\infty}}{{\text{lim}}}g_{{i}}z^{{i}}\cdot\underset{{j\rightarrow\infty}}{{\text{lim}}}h_{{j}}z^{{j}}=1$.

(\emph{ii}) The power double sequence of second type $\left\{ \left(P_{{A;g,h}}(z)\right)_{ij}\right\} $
is a regular method of summability.

\bigskip{}

\textbf{Proof}. The main structure of this proof is similar to the
one of the Theorem A. Let's show first that (\emph{ii}) implies (\emph{i}).
From (\emph{ii}) and the condition 3. of Silverman-Toeplitz theorem,
we have

\noindent $1={\underset{{i\rightarrow\infty}}{{\text{lim}}}\overset{{\infty}}{\underset{{j=0}}{\sum}}{a_{{ij}}g_{{i}}h_{{j}}z^{{i+j}}}}$.
This equals to $\underset{{i\rightarrow\infty}}{{\text{lim}}}g_{{i}}z^{{i}}\cdot\underset{{i\rightarrow\infty}}{{\text{lim}}}\overset{{\infty}}{\underset{{j=0}}{\sum}}{a_{{ij}}h_{{j}}z^{{j}}}=\underset{{i\rightarrow\infty}}{{\text{lim}}}g_{{i}}z^{{i}}\cdot\underset{{j\rightarrow\infty}}{{\text{lim}}}h_{{j}}z^{{j}}$,
where we also used that $\{a_{ij}\}$ as a linear operator preserves
limits by Silverman-Toeplitz theorem.

Now we show that (\emph{i}) implies (\emph{ii}). We will use again
Silverman-Toeplitz theorem and need to prove its three conditions:

1. ${\underset{{0\le i<\infty}}{{\text{sup}}}\overset{{\infty}}{\underset{{j=0}}{\sum}}{\left|\left(P_{{A;g,h}}(z)\right)_{ij}\right|}=\underset{{0\le i<\infty}}{{\text{sup}}}\overset{{\infty}}{\underset{{j=0}}{\sum}}{\left|a_{{ij}}g_{{i}}h_{{j}}z^{{i+j}}\right|}}$
by the definition of $P_{{A;g,h}}$.

\noindent The above supremum equals to
  \\ $\underset{{0\le i<\infty}}{{\text{sup}}}\left|g_{{i}}z^{{i}}\right|\cdot\overset{{\infty}}{\underset{{j=0}}{\sum}}{\left|a_{{ij}}h_{{j}}z^{{j}}\right|}\leq\underset{{0\le i<\infty}}{{\text{sup}}}\overset{{\infty}}{\underset{{j=0}}{\sum}}{|a_{{i,j}}|}\cdot\underset{{0\le i<\infty}}{{\text{sup}}}\left|g_{{i}}z^{{i}}\right|\cdot\underset{{0\le j<\infty}}{{\text{sup}}}{\left|h_{{j}}z^{{j}}\right|}<\infty$
  \\ since the first supremum is finite by the first condition of $\{a_{ij}\}$
being a regular method of summability and the other two are finite
by the existence of the limit in (\emph{i}).\bigskip{}

2. ${\underset{{i\rightarrow\infty}}{{\text{lim}}}\left(P_{{A;g,h}}(z)\right)_{ij}}={\underset{{i\rightarrow\infty}}{{\text{lim}}}a_{{ij}}g_{{i}}h_{{j}}z^{{i+j}}}$
for $j=0,1,2,\ldots$ by the definition.

\noindent The absolute value of the above limit is $\left|h_{{j}}z^{{j}}\cdot\underset{{i\rightarrow\infty}}{{\text{lim}}}a_{{ij}}g_{{i}}z^{{i}}\right|\leq\underset{{0\le i<\infty}}{{\text{sup}}}{|a_{{i,j}}|}\cdot\left|h_{{j}}z^{{j}}\cdot\underset{{i\rightarrow\infty}}{{\text{lim}}}g_{{i}}z^{{i}}\right|=0$
for $j=0,1,2,\ldots$

\noindent since the first limit is zero by the second condition of
$\{a_{ij}\}$ being a regular method of summability and the second
limit is finite by (\emph{i}).

\bigskip{}

3. ${\underset{{i\rightarrow\infty}}{{\text{lim}}}\overset{{\infty}}{\underset{{j=0}}{\sum}}{\left(P_{{A;g,h}}(z)\right)_{ij}}}={\underset{{i\rightarrow\infty}}{{\text{lim}}}\overset{{\infty}}{\underset{{j=0}}{\sum}}{a_{{ij}}g_{{i}}h_{{j}}z^{{i+j}}}}$
by the definition.

\noindent Starting with (\emph{i}) following the proof of necessary
condition backwards we have:

\noindent $1=\underset{{i\rightarrow\infty}}{{\text{lim}}}g_{{i}}z^{{i}}\cdot\underset{{j\rightarrow\infty}}{{\text{lim}}}h_{{j}}z^{{j}}=\underset{{i\rightarrow\infty}}{{\text{lim}}}g_{{i}}z^{{i}}\cdot\underset{{i\rightarrow\infty}}{{\text{lim}}}\overset{{\infty}}{\underset{{j=0}}{\sum}}{a_{{ij}}h_{{j}}z^{{j}}}={\underset{{i\rightarrow\infty}}{{\text{lim}}}\overset{{\infty}}{\underset{{j=0}}{\sum}}{a_{{ij}}g_{{i}}h_{{j}}z^{{i+j}}}}$.

And this finishes the proof that $\left\{ \left(P_{{A;g,h}}(z)\right)_{ij}\right\} $
is a regular method of summability by Silverman-Toeplitz theorem.

\bigskip{}

\textbf{Remarks}. (\emph{i}) implies $|z|=r_{g}=r_{h}$, where $r_{g}$
and $r_{h}$ are the radii of convergence. \hide{??? The requirement in
the Theorem B. that both sequences $\{g_{{i}}z^{{i}}\}$ and $\{h_{{i}}z^{{i}}\}$
must be convergent seems to be too restrictive but the condition (\emph{ii})
does not guarantee their convergence. There are examples of non-convergent
sequences $\{g_{{i}}z^{{i}}\}$ and $\{h_{{i}}z^{{i}}\}$ (??? for
all combinations of bounded and unbounded cases) and regular methods
of summability that map these sequences to convergent ones (= that
sum them). Then again by choosing $z=1$ one gets a counterexample
for each case.}

\bigskip{}

\textbf{Corollary}. From the proof it is clear that for $|z|<r_{g}$
and $|z|<r_{h}$, conditions 1. and 2. hold but the limit in 3. is
zero and we don't get a regular method of summability in that case.

\bigskip{}

\section{Boundedness Results}\mylabel{sec:6}

Assume ${A}\in{B_{{c}}}$ now, but not necessarily a regular method
of summability. It is clear that ${P_{{A;h}}^{{R}}(0)}\in{B_{{c}}}$
and ${P_{{A;h}}^{{R}}(1)}\in{B_{{c}}}$. Also ${A=P_{{A;g}}^{{C}}(1)\in{B_{{c}}}}$.
On the other hand, for a given $z$ we can ask: Does ${P_{{A;g}}^{{C}}(z)}\in{B_{{c}}}$
or ${P_{{A;h}}^{{R}}(z)}\in{B_{{c}}}$ hold? And because the conditions
1 and 2 clearly hold this is equivalent to: Does

\[
{\underset{{i\rightarrow\infty}}{{\text{lim}}}\overset{{\infty}}{\underset{{j=0}}{\sum}}{a_{{ij}}g_{i}z^{{i}}}}\textrm{ or }{\underset{{i\rightarrow\infty}}{{\text{lim}}}\overset{{\infty}}{\underset{{j=0}}{\sum}}{a_{{ij}}h_{j}z^{{j}}}}\textrm{exist?}\]

{ The next propositions provide an answer for these two kinds of
power matrices. For power double sequences of first type we have:}

\bigskip{}

{ \textbf{Proposition 2}(C). Let \textit{A} = \{ ${a_{{ij}}}$ \},
\textit{i} = 1, 2, \ldots{} , \textit{j} = 0, 1, 2, \ldots{} , be
a double sequence of complex scalars satisfying ${A}\in{B_{{c}}}$
If ${P_{{A;g}}^{{C}}(a)}\in{B_{{c}}}$, for some $a\neq0$, then ${P_{{A;g}}^{{C}}(z)}\in{B_{{c}}}$,
for all \textit{$z$} satisfying $\|z\| < \|a\|$.}

\textbf{Proof.} It is a straightforward verification of Silverman-Toeplitz
Theorem conditions.

\bigskip{}

{ \textbf{Proposition 2(R}). Let \textit{A} = \{ ${a_{{ij}}}$ \},
\textit{i} = 1, 2, \ldots{} , \textit{j} = 0, 1, 2, \ldots{} , be
a double sequence of complex scalars satisfying{ }${A}\in{B_{{c}}}$.
If ${P_{{A;h}}^{{R}}(a)}\in{B_{{c}}}$, for some $a\neq0$, then ${P_{{A;h}}^{{R}}(z)}\in{B_{{c}}}$,
for all \textit{$z$} satisfying $|z|<|a|$.}

{ \textbf{Proof.} From the above argument, we have ${P_{{A}}^{{R}}(0)}\in{B_{{c}}}$.
We only need to prove ${P_{{A}}^{{R}}(z)}\in{B_{{c}}}$, for all $|z|<|a|$
and $z\neq0$.}

{ From the hypothesis ${P_{{A}}^{{R}}(a)}\in{B_{{c}}}$, we have
} 
\begin{enumerate}
\item { ${\underset{{1\le i<\infty}}{{\text{sup}}}\overset{{\infty}}{\underset{{j=0}}{\sum}}{|a_{{ij}}a^{{j}}|}=M<\infty}$;} 
\item { ${\underset{{i\rightarrow\infty}}{{\text{lim}}}a_{{ij}}a^{{j}}}$
exists for \textit{j} = 1, 2, \ldots{} ;} 
\item { ${\underset{{i\rightarrow\infty}}{{\text{lim}}}\overset{{\infty}}{\underset{{j=0}}{\sum}}{a_{{ij}}a^{{j}}}}$
exists.} 
\end{enumerate}
{ We have to show} 
\begin{enumerate}
\item { ${\underset{{1\le i<\infty}}{{\text{sup}}}\overset{{\infty}}{\underset{{j=0}}{\sum}}{|a_{{ij}}z^{{j}}|}<\infty}$;} 
\item { ${\underset{{i\rightarrow\infty}}{{\text{lim}}}a_{{ij}}z^{{j}}}$
exists for \textit{j} = 1, 2, \ldots{} ;} 
\item { ${\underset{{i\rightarrow\infty}}{{\text{lim}}}\overset{{\infty}}{\underset{{j=0}}{\sum}}{a_{{ij}}z^{{j}}}}$
exists.} 
\end{enumerate}
{ In fact, from the first condition for ${P_{{A}}^{{R}}(a)}$ , we
obtain}

\[
{\underset{{1\le i<\infty}}{{\text{sup}}}\overset{{\infty}}{\underset{{j=0}}{\sum}}{|a_{{ij}}z^{{j}}|}=\underset{{1\le i<\infty}}{{\text{sup}}}\overset{{\infty}}{\underset{{j=0}}{\sum}}{|a_{{ij}}a^{{j}}|}\;|\frac{z}{a}|^{{j}}\le\underset{{1\le i<\infty}}{{\text{sup}}}\overset{{\infty}}{\underset{{j=0}}{\sum}}{|a_{{ij}}a^{{j}}|}=M<\infty}.\]

\bigskip{}

{ So ${P_{{A}}^{{R}}(z)}$ satisfies its first condition. Similarly,
from the Condition 2 of ${P_{{A}}^{{R}}(a)}$ , we have}

\[
{{\underset{{i\rightarrow\infty}}{{\text{lim}}}a_{{ij}}z^{{j}}}={\underset{{i\rightarrow\infty}}{{\text{lim}}}a_{{ij}}a^{{j}}\left(\frac{z}{a}\right)^{{j}}}=0\textrm{, for}j=1,2,\ldots.}\]

{ So ${P_{{A}}^{{R}}(z)}$ satisfies its second condition. Next we
show that ${P_{{A}}^{{R}}(z)}$ satisfies its condition 3 from Lemma
1.}

{ For any given ${\varepsilon>0}$, there exists \textit{N}, such
that ${|\frac{z}{a}|^{{N-1}}<\frac{\varepsilon}{4M}}$. From Conditions
1 and 2 above, there exists $K > 0$ such that for all \textit{m},
$n > K$, the following inequality holds} 

\[
{|a_{{mj}}a^{{j}}-a_{{nj}}a^{{j}}|<\frac{\varepsilon}{2N}}.\]

\noindent { Now for all $m,n>K$, we have} 

\begin{eqnarray*}
 &  & |\overset{{\infty}}{\underset{{j=0}}{\sum}}{a_{{mj}}z^{{j}}}-\overset{{\infty}}{\underset{{j=0}}{\sum}}{a_{{nj}}z^{{j}}}| \leqq\\
 & \leqq & |\overset{{N-1}}{\underset{{j=0}}{\sum}}{a_{{mj}}z^{{j}}}-\overset{{N-1}}{\underset{{j=0}}{\sum}}{a_{{nj}}z^{{j}}}|+|\overset{{\infty}}{\underset{{j=N}}{\sum}}{a_{{mj}}z^{{j}}}-\overset{{\infty}}{\underset{{j=N}}{\sum}}{a_{{nj}}z^{{j}}}|\\
 & \leqq & \overset{{N-1}}{\underset{{j=0}}{\sum}}{|a_{{mj}}a^{{j}}-a_{{nj}}a^{{j}}|}|\frac{z}{a}|^{{j}}+\\
 &  & |\frac{z}{a}|^{{N-1}}\left(\overset{{\infty}}{\underset{{j=N}}{\sum}}{|a_{{mj}}a^{{j}}|}|\frac{z}{a}|^{{j-N+1}}+\overset{{\infty}}{\underset{{j=N}}{\sum}}{|a_{{nj}}a^{{j}}|}|\frac{z}{a}|^{{j-N+1}}\right)\\
 & \leqq & \overset{{N-1}}{\underset{{j=0}}{\sum}}{|a_{{mj}}a^{{j}}-a_{{nj}}a^{{j}}|}+|\frac{z}{a}|^{{N-1}}\left(\overset{{\infty}}{\underset{{j=N}}{\sum}}{|a_{{mj}}a^{{j}}|}+\overset{{\infty}}{\underset{{j=N}}{\sum}}{|a_{{nj}}a^{{j}}|}\right)\\
 & < & \frac{N\varepsilon}{2N}+\frac{\varepsilon}{4M}\left(M+M\right)\\
 & = & \varepsilon.
\end{eqnarray*}

\noindent { This proposition is proved.}

\bigskip{}

{ Proposition 2(R) indicates that, the row power functions matrix
${P_{{A}}^{{R}}(z)}$ has a similar property to power series: If there
exists a number $a\neq0$, such that ${P_{{A}}^{{R}}(a)\in B_{{c}}}$,
then there exists a positive number ${r_{{A}}}$ such that, ${P_{{A}}^{{R}}(z)\in B_{{c}}}$,
for all ${|z|<r_{{A}}}$, and ${P_{{A}}^{{R}}(z)\notin B_{{c}}}$,
for all ${|z|>r_{{A}}}$. ${r_{{A}}}$ is called the radius of summability
of the matrix \textit{A}. The radius of summability of the matrix
\textit{A} is 0, if there does not exist a number $a\neq0$, such
that ${P_{{A}}^{{R}}(a)\in B_{{c}}}$; The radius of summability of
the matrix \textit{A} is $\infty$, if ${P_{{A}}^{{R}}(a)\in B_{{c}}}$
for all numbers \textit{a}.}

\bigskip{}

{ The following corollary follows immediately from Proposition 2(R)
and the above notations. }

{ \textbf{Corollary 3}. Let \textit{A} = \{ ${a_{{ij}}}$ \}, \textit{i}
= 1, 2, \ldots{} , \textit{j} = 0, 1, 2, \ldots{} , be a double
sequence of complex scalars. If ${A\in B_{{c}}}$ then ${r_{{A}}\geq1}$.
}

{ For any given row power matrix ${P_{{A}}^{{R}}(z)}$ , the entries
of any fixed row, \textit{i}, can be treated as the terms of a power
series }

\[
{\overset{{\infty}}{\underset{{j=0}}{\sum}}{a_{{ij}}z^{{j}}}}\]

{ Its radius of convergence is denoted by ${r_{{A}}^{{i}}}$ , for
\textit{i} = 1, 2, \ldots{} .}

\bigskip{}

{ \textbf{Proposition 4}. ${r_{{A}}\leqq\underset{{1\le i<\infty}}{{\text{inf}}}r_{{A}}^{{i}}}$
.}

{ Proof. For any given ${|z|<r_{{A}}}$ , ${P_{{A}}^{{R}}(z)\in B_{{c}}}$,
we have ${\underset{{i\rightarrow\infty}}{{\text{lim}}}\overset{{\infty}}{\underset{{j=0}}{\sum}}{a_{{ij}}z^{{j}}}}$
exists. The series ${\overset{{\infty}}{\underset{{j=0}}{\sum}}{a_{{ij}}z^{{j}}}}$
is convergent, for $i=1,2,...$ . It implies that ${|z|\leqq r_{{A}}^{{i}}}$
, for \textit{i} = 1, 2, \ldots{} . It completes the proof of this
proposition.}

\bigskip{}

For power double sequences of second type we have:

\global\long\global\long\global\long\def\labelenumi{(\alph{enumi})}
 \global\long\global\long\global\long\def\labelenumii{\alph{enumii}}
 \global\long\global\long\global\long\def\labelenumiii{\alph{enumiii}}

\bigskip{}

{ \textbf{Proposition} 2(II). Let $\{a_{ij}\}$, $i=0,1,2,\ldots$,
$j=0,1,2,\ldots$, be a double sequence of complex scalars and $g$
a complex power series. }

{ If ${P_{{A;g}}(a)}\in{B_{{c}}}$, for some $a\neq0$, then ${P_{{A;g}}(z)}\in{B_{{c}}}$,
for all \textit{$z$} satisfying $|z|<|a|$.}

{ \textbf{Proof}. It is clear that ${P_{{A;g}}(0)}\in{B_{{c}}}$
and we only need to prove ${P_{{A;g}}(z)}\in{B_{{c}}}$, for all $0<|z|<|a|$.}

{ From the hypothesis ${P_{{A;g}}(a)}\in{B_{{c}}}$, we have } 
\begin{enumerate}
\item { ${\underset{{0\le i<\infty}}{{\text{sup}}}\overset{{\infty}}{\underset{{j=0}}{\sum}}{|a_{{ij}}g_{{i+j}}a^{{i+j}}|}=M<\infty}$;} 
\item { ${\underset{{i\rightarrow\infty}}{{\text{lim}}}a_{{ij}}g_{{i+j}}a^{{i+j}}}$
exists for \textit{j} = 1, 2, \ldots{} ;} 
\item { ${\underset{{i\rightarrow\infty}}{{\text{lim}}}\overset{{\infty}}{\underset{{j=0}}{\sum}}{a_{{ij}}g_{{i+j}}a^{{i+j}}}}$
exists.} 
\end{enumerate}
{ We have to show that the above three conditions are also true for
${P_{{A;g}}(z)}$:} 
\begin{enumerate}
\item { ${\underset{{0\le i<\infty}}{{\text{sup}}}\overset{{\infty}}{\underset{{j=0}}{\sum}}{|a_{{ij}}g_{{i+j}}z^{{i+j}}|}<\infty}$;} 
\item { ${\underset{{i\rightarrow\infty}}{{\text{lim}}}a_{{ij}}g_{{i+j}}z^{{i+j}}}$
exists for \textit{j} = 0,1,2, \ldots{} ;} 
\item { ${\underset{{i\rightarrow\infty}}{{\textrm{lim}}}\overset{{\infty}}{\underset{{j=0}}{\sum}}{a_{{ij}}g_{{i+j}}z^{{i+j}}}}$
exists.} 
\end{enumerate}
{ In fact, from the condition (\emph{a}) for ${P_{{A;g}}(a)}$ ,
we obtain}

\begin{eqnarray*}
\underset{{0\le i<\infty}}{{\textrm{sup}}}\overset{{\infty}}{\underset{{j=0}}{\sum}}{|a_{{ij}}g_{{i+j}}z^{{i+j}}|} & = & \underset{{0\le i<\infty}}{\textrm{sup}}\overset{{\infty}}{\underset{{j=0}}{\sum}}{|a_{{ij}}g_{{i+j}}a^{{i+j}}|}|\frac{z}{a}|^{{i+j}}\\
         & \le & \underset{{0\le i<\infty}}{\textrm{sup}}\overset{{\infty}}{\underset{{j=0}}{\sum}}{|a_{{ij}}g_{{i+j}}a^{{i+j}}|} = M < \infty
\end{eqnarray*}

{ So ${P_{{A;g}}(z)}$ satisfies its first condition. Similarly,
from the condition (\emph{b}) of ${P_{{A;g}}(a)}$ , we have}

\[
{\underset{{i\rightarrow\infty}}{{\textrm{lim}}}a_{{ij}}g_{{i+j}}z^{{i+j}}}={\underset{{i\rightarrow\infty}}{{\textrm{lim}}}a_{{ij}}g_{{i+j}}a^{{i+j}}\left(\frac{z}{a}\right)^{{i+j}}}=0\textrm{, for}j=0,1,2,\ldots\]

{ So ${P_{{A;g}}(z)}$ satisfies its second condition. Next we show
that ${P_{{A;g}}(z)}$ satisfies its condition (\emph{c}).}

{ For any given ${\varepsilon>0}$, there exists \textit{N}, such
that ${|\frac{z}{a}|^{{N-1}}<\frac{\varepsilon}{4M}}$. From conditions
(\emph{a}) and (\emph{b}) above, where we set $k=i+j$, there exists
\textit{$K>0$} such that ${|\frac{z}{a}|^{{K}}<\frac{\varepsilon}{4MN}}$
and for all \textit{m}, $n > K$, the following inequalities
hold (without loss of generality $m<n$) }

\[
{\overset{{\infty}}{\underset{{j=0}}{\sum}}{|a_{{m-j,j}}g_{{m}}a^{{m}}|}<M},\]

\[
{\overset{{\infty}}{\underset{{j=0}}{\sum}}{|a_{{n-j,j}}g_{{n}}a^{{n}}|}<M},\]

\[
{|a_{{m-j,j}}g_{{m}}a^{{m}}-a_{{n-j,j}}g_{{n}}a^{{n}}|<\frac{\varepsilon}{4N}}.\]

\noindent { Now for all \textit{m},  $n > K$, we have}

\begin{eqnarray*}
 &  & \left|\overset{{\infty}}{\underset{{j=0}}{\sum}}{a_{{m-j,j}}g_{{m}}z^{{m}}}-\overset{{\infty}}{\underset{{j=0}}{\sum}}{a_{{n-j,j}}g_{{n}}z^{{n}}}\right| \leqq\\
 & \leqq & \left|\overset{{N-1}}{\underset{{j=0}}{\sum}}{a_{{m-j,j}}g_{{m}}z^{{m}}}-\overset{{N-1}}{\underset{{j=0}}{\sum}}{a_{{n-j,j}}g_{{n}}z^{{n}}}\right|+\left|\overset{{\infty}}{\underset{{j=N}}{\sum}}{a_{{m-j,j}}g_{{m}}z^{{m}}}-\overset{{\infty}}{\underset{{j=N}}{\sum}}{a_{{n-j,j}}g_{{n}}z^{{n}}}\right|\\
 & \leqq & \overset{{N-1}}{\underset{{j=0}}{\sum}}{\left|a_{{m-j,j}}g_{{m}}a^{{m}}|\frac{z}{a}|^{{m}}-a_{{n-j,j}}g_{{n}}a^{{n}}|\frac{z}{a}|^{{n}}\right|}+\\
 &  & |\frac{z}{a}|^{{N-1}}\left(\overset{{\infty}}{\underset{{j=N}}{\sum}}{\left|a_{{m-j,j}}g_{{m}}a^{{m}}\right|}|\frac{z}{a}|^{{m-N+1}}+\overset{{\infty}}{\underset{{j=N}}{\sum}}{\left|a_{{n-j,j}}g_{{n}}a^{{n}}\right|}|\frac{z}{a}|^{{n-N+1}}\right)\\
 & \leqq & \overset{{N-1}}{\underset{{j=0}}{\sum}}{\left(\left|a_{{m-j,j}}g_{{m}}a^{{m}}-a_{{n-j,j}}g_{{n}}a^{{n}}\right|\cdot|\frac{z}{a}|^{{m}}+\left|a_{{n-j,j}}g_{{n}}a^{{n}}\right|\cdot(1-|\frac{z}{a}|^{{n-m}})\cdot|\frac{z}{a}|^{{m}}\right)}+\\
 &  & |\frac{z}{a}|^{{N-1}}\left(\overset{{\infty}}{\underset{{j=N}}{\sum}}{\left|a_{{m-j,j}}g_{{m}}a^{{m}}\right|}+\overset{{\infty}}{\underset{{j=N}}{\sum}}{\left|a_{{n-j,j}}g_{{n}}a^{{n}}\right|}\right)\\
 & < & \frac{N\varepsilon}{4N}+NM\frac{\varepsilon}{4MN}+\frac{\varepsilon}{4M}\left(M+M\right)\\
 & = & \varepsilon.\end{eqnarray*}

\noindent { This proposition is proved.}

\bigskip{}
For power double sequences of third type we have:

\bigskip{}

{ \textbf{Proposition} 2(III). Let $\{a_{ij}\}$, $i=0,1,2,\ldots$,
$j=0,1,2,\ldots$, be a double sequence of complex scalars and $g,h$
complex power series. If ${P_{{A;g,h}}(a)}\in{B_{{c}}}$, for some
$a\neq0$, then ${P_{{A;g,h}}(z)}\in{B_{{c}}}$, for all \textit{$z$}
satisfying $|z|<|a|$.}

{ \textbf{Proof}. This proof is very similar to the one of the Proposition
2(II) above and will be omitted.\bigskip{}

\textbf{Remarks}. Propositions 2(II) and 2(III) might seem to generalize
the Proposition 2(C) or Proposition 2(R) but they don't. We cannot
write ${P_{{A}}^{{R}}(z)}$ nor ${P_{{A}}^{{C}}(z)}$ in the form
of $P_{{A;g}}$ or $P_{{A;g,h}}$. But we can do so for ${P_{{A}}(z)}$
by choosing $g_{i}=1$ or $g_{i}=h_{i}=1$ respectively.

Obviously ${P_{{A;g}}(0)}\in{B_{{c}}}$, ${A=P_{{A;g}}(1)}\in{B_{{c}}}$,
${P_{{A;g,h}}(0)}\in{B_{{c}}}$, and ${A=P_{{A;g,h}}(1)}\in{B_{{c}}}$.
Then for these power matrices ${r_{{A}}\geq1}$.

\bigskip{}

\section{Examples}\mylabel{sec:7}

For the following examples, one can check the conditions listed.

\noindent { \textbf{Example 1}. Let \textit{A} = \{\textit{i} + \textit{j}\},
\textit{i} = 1, 2, \ldots{}, \textit{j} = 0, 1, 2, \ldots{} :} 
\begin{enumerate}
\item { ${A\notin B_{{c}}}$; } 
\item { ${r_{{A}}^{{i}}=1}$, for \textit{i} = 1, 2, \ldots{} ;} 
\item { ${r_{{A}}=1}$. } 
\end{enumerate}
{ \textbf{Example 2}. Let \textit{A} = \{ ${\frac{1}{i+j}}$ \},
\textit{i} = 1, 2, \ldots{} , \textit{j} = 0, 1, 2, \ldots{} :} 
\begin{enumerate}
\item { ${A\notin B_{{c}}}$; } 
\item { ${r_{{A}}^{{i}}=1}$, for \textit{i} = 1, 2, \ldots{} ;} 
\item { ${r_{{A}}=1}$. } 
\end{enumerate}
{ \textbf{Example 3}. Let \textit{A} = \{ ${(i+j)!}$ \}, \textit{i}
= 1, 2, \ldots{} , \textit{j} = 0, 1, 2, \ldots{} :} 
\begin{enumerate}
\item { ${A\notin B_{{c}}}$; } 
\item { ${r_{{A}}^{{i}}=0}$, for \textit{i} = 1, 2, \ldots{} ;} 
\item { ${r_{{A}}=0}$. } 
\end{enumerate}
{ \textbf{Example 4}. Let \textit{A} = \{ ${\frac{1}{(i+j)!}}$ \},
\textit{i} = 1, 2, \ldots{} , \textit{j} = 0, 1, 2, \ldots{} :} 
\begin{enumerate}
\item { ${A\in B_{{c}}}$, with $M={\overset{{\infty}}{\underset{{j=\;0}}{\sum}}{\frac{1}{j!}}}$;
} 
\item { ${r_{{A}}^{{i}}=\infty}$, for \textit{i} = 1, 2, \ldots{}; } 
\item { ${r_{{A}}=\infty}$. } 
\end{enumerate}
{ \textbf{Example 5}. For a given $b > 0$, let \textit{A}
= \{ ${b^{{i+j}}}$ \}, \textit{i} = 1, 2, \ldots{} , \textit{j}
= 0, 1, 2, \ldots{} :} 
\begin{enumerate}
\item { ${A\in B_{{c}}}$, if \textit{b} < 1, with $M={\frac{1}{1-b}}$
, and \textit{A} $\notin{B_{{c}}}$ , if \textit{b} $\geq$ 1; \ } 
\item { ${r_{{A}}^{{i}}}={b^{{-1}}}$ , for \textit{i} = 1, 2, \ldots{}
;} 
\item { ${r_{{A}}}={b^{{-1}}}$ . } 
\end{enumerate}
{ \textbf{Example 6}. Let ${[p_{{k}}]}$ be a sequence of positive
numbers satisfying ${\overset{{\infty}}{\underset{{k=1}}{\sum}}{p_{{k}}}}<\infty$,
and let ${P_{{i}}=\overset{{i}}{\underset{{k=1}}{\sum}}{p_{{k}}}}$.
Define ${A=a_{{ij}},\: i=1,2,...,\: j=0,1,2,...}$, as follows}

\[
{{a_{ij}=\left\{ \begin{matrix}\frac{1}{P_{i}},j<i\\
0,j\ge i\end{matrix}\right.}{:}}\]

\begin{enumerate}
\item { ${A\in B_{{c}}}$, with $M={\frac{1}{p_{{1}}}}$ ; \ } 
\item { ${r_{{A}}^{{i}}=\infty}$, for \textit{i} = 1, 2, \ldots{};} 
\item { ${r_{{A}}=1}$;} 
\item { \textit{A} is not a regular method of summability.} 
\end{enumerate}
{ \textbf{Example 7}. Let ${[p_{{k}}]}$ be a sequence of positive
numbers satisfying ${\overset{{\infty}}{\underset{{k=1}}{\sum}}{p_{{k}}}}=\infty$,
and let ${P_{{i}}=\overset{{i}}{\underset{{k=1}}{\sum}}{p_{{k}}}}$.
Define $A=a_{{ij}}$, as in Example 6. Then:} 
\begin{enumerate}
\item { ${A\in B_{{c}}}$, with $M={\frac{1}{p_{{1}}}}$ ; \ } 
\item { ${r_{{A}}^{{i}}=\infty}$, for \textit{i} = 1, 2, \ldots{};} 
\item { ${r_{{A}}}={\text{sup}\{b>0:\underset{{i\rightarrow\infty}}{{\text{lim}}}\frac{b^{{i}}}{P_{{i}}}}$
exists\} $<\infty$;} 
\item { \textit{A} is a regular method of summability.} \bigskip{}

\end{enumerate}
{ \textbf{Proposition 6}. Let $f(z)$ be an analytic function with
power series expansion ${\overset{{\infty}}{\underset{{k=0}}{\sum}}{p_{{k}}z^{{k}}}}$
about point 0 with the radius of convergence \textit{r}. Then, for
any $s=[s_{{1}},\; s_{{2}},\; s_{{3}},\;\ldots]\in c$ and for any
$|z|<r$, the sequence $[\overset{{n}}{\underset{{k=0}}{\sum}}{p_{{k}}s_{{k+1}}z^{{k}}}]$,
\textit{n} = 1, 2, \ldots{} , is also convergent.}

{ Proof. Define \textit{A} = \{ ${a_{{ij}}}$ \}, \textit{i} = 1,
2, 3, \ldots{} , \textit{j} = 0, 1, 2, \ldots{} , as follows:}

\[
{a_{ij}=\left\{ \begin{matrix}p_{j},j<i\\
0,j\ge i\end{matrix}\right.}\]

{ One can check that \textit{A}, \textit{A}(\textit{z}) have the
following properties:} 
\begin{enumerate}
\item { ${A\in B_{{c}}}$, with ${M=\overset{{\infty}}{\underset{{k=0}}{\sum}}{p_{{k}}}}$;
} 
\item { ${r_{{A}}^{{i}}=\infty}$, for \textit{i} = 1, 2, \ldots{};} 
\item { ${r_{{A}}=r}$;} 
\item { ${A(z)\in B_{{c}}}$ .} 
\end{enumerate}
{ Then this proposition follows immediately from the property (d).}




\bibliographystyle{plain}

\bigskip
\noindent
\parbox[t]{.48\textwidth}{
Jinlu Li \\
Department of Mathematics \\
Shawnee State University \\
940 Second Street \\
Portsmouth, OH 45662, USA \\
jli@shawnee.edu 
} \hfill
\parbox[t]{.48\textwidth}{
Robert Mendris \\
Department of Mathematics \\
Shawnee State University \\
940 Second Street \\
Portsmouth, OH 45662, USA \\
rmendris@shawnee.edu 
} \hfill

\end{document}